\documentclass[a4paper, 10pt]{article}
\usepackage{latexsym,amsthm,amssymb,amscd}
\usepackage{color,graphicx}
\usepackage{framed}
\author{{{\bf  S. Golalizadeh and }}  {\bf N. Soltankhah \thanks{Corresponding author: soltan@alzahra.ac.ir, soltankhah.n@gmail.com }} \\
{\footnotesize  {\bf Department of Mathematical Sciences}}\\  {\footnotesize {\bf Alzahra University}}\\ {\footnotesize {\bf P.O. Box 19834}}\\ {\footnotesize {\bf Tehran, I.R. Iran.}}}
\title {The minimum possible volume size of  $\mu$-way $(v,k,t)$ trades  }
\date{}
\newtheorem{theorem}{Theorem}

\newtheorem{lemma}{Lemma}
\newtheorem{corollary}{Corollary}
\newtheorem{example}{Example}

\begin{document}
\maketitle
\vspace*{-9cm}
\begin{framed}
\textbf{PLEASE CITE AS:}
  S. Golalizadeh and  N. Soltankhah. The minimum possible volume size of  $\mu$-way $(v,k,t)$ trades. Utilitas Mathematica, 111 (2019), 211-224.
\end{framed}
\vspace*{7cm}
\begin{abstract}
A $\mu$-way $(v, k, t)$ trade is a pair $T=(X,\{T_1,T_2,\cdots, T_{\mu}\})$  such that for each $t$-subset of $v$-set $X$ the number of blocks containing this $t$-subset is the same in each $T_i$ $(1\leq i\leq \mu)$. In the other words for  each $1\leq i< j\leq \mu$,  $(X,\{T_i,T_j\})$ is a $(v,k,t)$ trade. There are many questions concerning $\mu$-way trades. The main question is about the minimum volume and minimum foundation size of $\mu$-way $(v,k,t)$ trades. In this paper, we determine the minimum volume and minimum foundation size of  $\mu$-way $(v,t+1,t)$ trades for each  integer number $\mu\geq 3$ and $t=2$.
\end{abstract}

\section{Introduction and preliminaries}
A $(v,k,t)$ trade is a pair $T=(X,\{T_1,T_2\})$ in which $T_1$ and $T_2$ are two disjoint collections of $k$-subsets (called blocks) of $X$ such that for every $t$-subset $A$, the number of occurences of $A$ in $T_1$ is the same as the number of occurences of $A$ in $T_2$. We will use the notation $(X;T_1,T_2)$ instead of $T=(X,\{T_1,T_2\})$.
A $(v, k, t)$ trade is called $(v, k, t)$ Steiner trade if any
$t$-subset of $X$ occurs in at most once in $T_1(T_2)$. A $(v,k,t)$ trade is called $d$-homogeneous if the number of occurences of each element of $X$ in $T_1$ ($T_2$) is exactly $d$. In a $(v, k, t)$ trade, both collections of blocks
must cover the same set of elements. This set of elements is called the foundation
of the trade and is denoted by found(T). Also, the number of blocks in $T_1$ is the same as the number of blocks in $T_2$ which is called the volume of the trade.\\

Trades are  useful combinatorial objects with many  applications in various areas of combinatorial design theory. The concept of trade was first introduced ($1960$) in paper \cite{hedayat}. Later,  Steiner trades  are used and renamed by Milici and Quattrocchi  ($1986$) with the name of DMB (disjoint and mutually balanced). 
 Recently, a generalization for the concept of trade has been introduced in \cite{rashidi} by the name of $\mu$-way trades as follows:

A $\mu$-way $(v, k, t)$ trade is a pair $T=(X,\{T_1,T_2,\cdots, T_{\mu}\})$ such that for each $i\neq j$, $(X;T_i,T_j)$ is a $(v,k,t)$ trade. \\

There exist many questions concerning $\mu$-way trades. Some of the most important questions are about the minimum volume and minimum foundation size and the set of all possible volume sizes of  $\mu$-way trades. Not much is known for the mentioned questiones about $\mu$-way $(v,k,t)$ trades for $\mu\geq 3$ and most of the papers have been focused on the case $\mu=2$, see  \cite{asgari,hedayat2,hwang,mahmood}. $\mu$-way trade is defined as different name as $N$-legged trade, before, see \cite{forbes}.  Some questions have been   answered   about the existence and non-existence of $3$-way $(v,k,t)$ trades for some special values of $k$ and $t$ in \cite{rashidi}. Also, some results on the  $3$-way $d$-homogeneous $(v,3,2)$ Steiner trades for some special values of $d$ have been obtained   in \cite{amjadi}. \\ 

In this paper,  we determine the minimum  volume and the minimum foundation size of  $\mu$-way $(v,t+1,t)$ trades for each integer number $\mu\geq 3$ and $t=2$. This problem has been solved by Hwang \cite{hwang} for $\mu=2$ and each integer number $t$ by an easy induction on $t$. But to obtain a similar result about $\mu$-way $(v,3,2)$ trades for each $\mu\geq 3$, we need more information about the structure of $\mu$-way trades of minimum volume.\\

In this  paper, we concern with  $(2,3)$-packing designs. These objects are very useful in the study of $(v,3,2)$ trades and their generalization. Conversely, $(v,3,2)$ trades and their generalizations can be used in the study of the intersection problem of packings and their large sets. \\

Let $V$ be a set of $v$ points. A $(2,3)$-packing  on $V$ is a pair $(V,\mathcal{A})$, where $\mathcal{A}$ is  a set of $3$-subsets (called blocks) of $V$, such that each pair of distinct points from $V$ appears together in at most one block. The leave of a $(2,3)$-packing $(V,\mathcal{A})$   is the graph
$(V, E)$ where $E$ consists of all the pairs which do not appear in any block of $\mathcal{A}$.\\

A $(2,3)$-packing $(V,\mathcal{A})$ is called maximum if there does not exist any $(2,3)$-packing $(V,\mathcal{B})$ with $|\mathcal{A}|<|\mathcal{B}|$. The maximum packings of the same order all have one things in common: the leave. If $(V,\mathcal{A})$ is a maximum $(2,3)$-packing of order $v$, according to \cite{lindner} the leave is\\
$(1)$ a $1$-factor if $v\equiv 0,2$ (mod $6$),\\
$(2)$ a $4$-cycle if $v\equiv 5$ (mod $6$),\\
$(3)$ a tripole, that is a spanning graph with each vertex having odd degree and containing $\frac{(v+2)}{2}$ edges if $v\equiv 4$ (mod $6$) and\\
$(4)$ the empty set if $v\equiv 1,3$ (mod $6$). In this case a $(2,3)$-packing is called Steiner triple system of order $v$ and the number of blocks is denoted by $t_v$. It has been shown that $t_v=\frac{v(v-1)}{6}$.\\

 Two $(2,3)$-packings $(V,\mathcal{A})$ and $(V,\mathcal{B})$ are called disjoint and mutually balanced (DMB)  if\\
$(1)$ $\mathcal{A}\cap \mathcal{B}=\emptyset$,\\
$(2)$ any given pair of distinct elements of $V$ is contained in a block of $\mathcal{A}$ if and only if it is
contained in a block of $\mathcal{B}$. On the other words two $(2,3)$-packings $(V,\mathcal{A})$ and $(V,\mathcal{B})$ have the same leave.\\
 A set of more than two $(2,3)$-packings is called DMB if each pair of them is DMB. A set of $m$ $(2,3)$-packings $(V,\mathcal{A}_i)$ for $i=1,2,\cdots,m$  intersect in $k$ blocks if $\mathcal{A}_i\cap \mathcal{A}_j=A$ for each $i\neq j$ and $|A|=k$.  By $\varphi(v)$,  we denote the maximum number of disjoint mutually balanced (DMB) $(2,3)$-packings of order $v$. Also  $D(v,k)$ will be denote the maximum number of maximum $(2,3)$-packings on $v$ points with the same leave such that any two of them have exactly $k$ blocks in common, these $k$ blocks are contained in each of the maximum $(2,3)$-packings.  The following upper bound on $\varphi(v)$  has been proved in \cite{schellenberg}.
\begin{theorem}\label{0}\cite{schellenberg}
$\varphi(v)\leq v-2$ for $v\equiv 1,3$ (mod $6$); $\varphi(v)\leq v-4$ for $v\equiv 0,2,5$ (mod $6$); and $\varphi(v)\leq v-6$ for $v\equiv 4$ (mod $6$). Further, except $v\equiv 4$ (mod $6$), the upper bound is acceived only if the packings are maximum.
\end{theorem}
\begin{theorem}\label{1}
\begin{itemize}
\item[1)] For $v\equiv 1,3$ (mod $6$) and $v\neq 7$, $\varphi(v)=v-2$. Also $\varphi(7)=2$\cite{lu1,lu2,teirlink}.
\item[2)] For $v\equiv 0,2$ (mod $6$), $\varphi(v)=v-4$\cite{chen,chen4,lei1}.
\item[3)] For $v\equiv 4$ (mod $6$), $\varphi(v)=v-6$\cite{cao2, cao3, cao4, chang, chen, ji}.
\item[4)] For $v\equiv 5$ (mod $6$), $\varphi(v)=v-4$\cite{cao1,chen2, chen3,schellenberg}.
\end{itemize}
\end{theorem}
Let $(X, \{T_1,T_2, \cdots,  T_{\mu}\})$ be a $\mu$-way
$(v, k, t)$ trade of volume $m$, and $x,y\in$ found($T$).
Let us denote by $r_x$ and $\lambda_{xy}$ the number of blocks in $T_i$ $(1\leq i \leq \mu)$ which contains the element $x$ and the number of blocks contains the pair $x,y$, respectively.  The set of blocks in $T_i$ $(1\leq i \leq \mu)$ which contains
$x \in$ found(T) is denoted by $T_{ix}$ $(1\leq i \leq \mu)$ and the set of remaining blocks by $T'_{ix}$ $(1\leq i \leq \mu)$.\\

As an analogue of a lemma in \cite{hwang}, we have the following lemma for $\mu$-way trades.
\begin{lemma}\label{32}
Let $T=(X,\{T_1,T_2,\cdots, T_{\mu}\})$ be a $\mu$-way $(v,k,t)$ trade of volume $m$ and $x\in $found(T) with $r_x<m$,
\begin{enumerate}
\item
$T_x = (X,\{T_{1x},T_{2x}, \cdots , T_{\mu x}\})$ is a $\mu$-way
$(v,k, t-1)$ trade of volume $r_x$,
\item
$T'_x = (X\setminus \{x\},\{T'_{1x},\cdots, T'_{\mu x}\})$ is a $\mu$-way $(v-1, k, t-1)$
trade of volume $m-r_x$.
\end{enumerate}
\end{lemma}
In the throughout of this paper, the notation $xA$ will denote the set that its elements are the union of the set $\{x\}$ with each element of $A$.
\section{Upper bounds on the numbers $\mathcal{D}(v,k)$}
In this section we obtain  upper bounds on the numbers $\mathcal{D}(v,k)$ for some special values of $k$ that will be useful to obtain main result. Note that, if there exists a set of $m$ $(2,3)$-packings on $v$-set  $V$ with the same leave and $b_v$ blocks which intersect in $k$ blocks,  the remaining different blocks form  a $m$-way $(v,3,2)$ Steiner trade of volume $b_v-k$. Conversely, If $(X,\{T_1,T_2,\cdots, T_m\})$ is a $m$-way $(v,3,2)$ Steiner trade, then the set $\{T_1,T_2,\cdots, T_m\}$ forms a set of $m$ DMB $(2,3)$-packings.
\begin{example}\label{pp}
Let $X=\{1,2,3,4,5,6,7\}$ and $\mathcal{A}_1$, $\mathcal{A}_2$ and $\mathcal{A}_3$ be the following sets.
$$\mathcal{A}_1=\{123,145,167,246,257,347,356\},$$
$$\mathcal{A}_2=\{123,146,157,247,256,345,367\},$$
$$\mathcal{A}_3=\{123,147,156,245,267,346,357\}.$$
It is easy to check that $(X,\mathcal{A}_1)$, $(X,\mathcal{A}_2)$ and $(X,\mathcal{A}_3)$ are three STS$(7)$s which  intersect in one block. Then $(X,\{T_1,T_2,T_3\})$ is a $3$-way $(7,3,2)$ Steiner trade of volume $6=7-1$ in which $T_i=\mathcal{A}_i\setminus \{123\}$ for $i=1,2,3$ (also the set $\{T_1,T_2,T_3\}$ forms a set of $3$ DMB $(2,3)$-packings).
\end{example}
\begin{example}\label{ppp}
Let $X=\{1,2,3,4,5,6,7,8,9\}$ and $\mathcal{A}_1$, $\mathcal{A}_2$ and $\mathcal{A}_3$ be the following sets.
$$\mathcal{A}_1=\{123,145,167,189,246,258,279,349,357,368,478,569\},$$
$$\mathcal{A}_2=\{123,146,158,179,247,259,268,345,369,378,489,567\},$$
$$\mathcal{A}_3=\{123,147,159,168,245,267,289,348,356,379,469,578\}.$$
It is readily checked that $(X,\mathcal{A}_1)$, $(X,\mathcal{A}_2)$ and $(X,\mathcal{A}_3)$ are three STS$(9)$s which  intersect in one block. Then $(X,\{T_1,T_2,T_3\})$ is a $3$-way $(9,3,2)$ Steiner trade of volume $11=12-1$ in which $T_i=\mathcal{A}_i\setminus \{123\}$ for $i=1,2,3$ (also the set $\{T_1,T_2,T_3\}$ forms a set of $3$ DMB $(2,3)$-packings).
\end{example}
\begin{example}\label{pppp}
Let $X=\{1,2,3,4,5,6,7,8,9\}$ and $\mathcal{A}_1$, $\mathcal{A}_2$, $\mathcal{A}_3$ and $\mathcal{A}_4$ be the following sets.
$$\mathcal{A}_1=\{123,145,167,189,246,258,279,349,357,368,478,569\},$$
$$\mathcal{A}_2=\{123,145,167,189,259,247,268,356,348,379,469,578\},$$
$$\mathcal{A}_3=\{123,145,167,189,249,256,278,347,358,369,468,579\},$$
$$\mathcal{A}_4=\{123,145,167,189,248,257,269,346,359,378,479,568\}.$$
It is readily checked that $(X,\mathcal{A}_1)$, $(X,\mathcal{A}_2)$, $(X,\mathcal{A}_3)$ and $(X,\mathcal{A}_4)$ are four STS$(9)$s which intersect in four blocks. Then $(X\setminus \{1\},\{T_1,T_2,T_3,T_4\})$ is a $4$-way $(8,3,2)$ Steiner trade of volume $8=12-4$ in which $T_i=\mathcal{A}_i\setminus \{123,145,167,189\}$ for $i=1,2,3,4$ (also the set $\{T_1,T_2,T_3,T_4\}$ forms a set of $4$ DMB $(2,3)$-packings).
\end{example}
\begin{theorem}\label{4}
\begin{itemize}
\item[i.]
If $v\equiv 1$ or $3$ (mod $6$), and $v\geq 13$, then $\mathcal{D}(v,1)\leq v-5$, $\mathcal{D}(7,1)=3$  and $\mathcal{D}(9,1)=3$.
\item[ii.]
If $v\equiv 0,2$ (mod $6$) and $v\geq 8$, $\mathcal{D}(v,1)\leq v-6$.
\end{itemize}
\begin{proof}
Let  $X$ be a set on $v$ points. Suppose that $(X,\mathcal{A}_1)$, $(X,\mathcal{A}_2)$, $\cdots$ and $(X,\mathcal{A}_t)$ are $t$ Steiner triple systems with exactly one block in common. Without loss of generality, consider  the block $\{1,2,3\}$ as the only common block. Assume that $w$ is an element of $X\setminus \{1,2,3\}$. The third element in block containing the pair $1w$ must be distinct in each of $t$ systems and this element  belongs to $X\setminus \{1,2,3,w\}$. Then $t\leq v-4$.  Suppose that $t=v-4$. Let $F^j_i$ be the $1$-factor on the set $X\setminus \{j\}$ such that $jF^j_i\subset \mathcal{A}_i$ for $i=1,2,\cdots,v-4$. Let $C_j=\bigcup_{k=1,k\neq j}^v F^k_{v-4}$. It is easy to see that the number of distinct elements of $C_j$ to form $ab$ where $a,b\notin\{1,2,3,j\}$ is at least  $(v-1)(\frac{v-1}{2}-4)+6$ and the number of distinct elements of $\bigcup_{i=1}^{v-5} F^j_i$ to the mentioned form is $(v-5)(\frac{v-1}{2}-3)$. Since for $v\geq 9$, $(v-1)(\frac{v-1}{2}-4)+6>(v-5)(\frac{v-1}{2}-3)$, then there exists at least an element for example $xy$ such that $xy\in C_j\setminus \bigcup_{i=1}^{v-5} F^j_i$. Since the pairs to form $xa$ such that $jxa\in \mathcal{A}_i$ for $i=1,2,\cdots,v-5$ are distinct and $jxy\notin \mathcal{A}_{v-4}$, then there does not exist any element for block containing the pair $jx$ in the Steiner triple system $(X,\mathcal{A}_{v-4})$ and this is a contradiction. \\

Suppose that $X$ is a set on $v$ points where $v\equiv 0,2$ (mod $6$). Let $(X,\mathcal{A}_1)$, $(X,\mathcal{A}_2)$, $\cdots$, $(X,\mathcal{A}_t)$ be the maximum packings with the $1$-factor $F$ as the common leave. Suppose that $123$ is the only common block of these packings. We choice an element $w\in X\setminus \{1,2,3\}$ and $1w\notin F$. Assume that $x$ and $y$ are the elements of the set $X$ such that  $1x, wy\in F$. The third element in block containing the pair $1w$ must be distinct in each of $t$ packings and this element  belongs to $X\setminus \{1, 2, 3, w, x, y\}$. Then $t\leq v-6$.
\end{proof}
\end{theorem}
\begin{theorem}\label{6}
\begin{itemize}
\item[i.]
If $v\equiv 1,3$ (mod $6$), then $\mathcal{D}(v,k)\leq v-6$, for $k\geq 2$ and $k\neq \frac{v-1}{2}$.
\item[ii.]
If $v\equiv 0,2$ (mod $6$), then $\mathcal{D}(v,k)\leq v-8$ for $k\geq 2$ and  $k\neq \frac{v-2}{2}$.
\end{itemize}
\begin{proof}
We prove the Case $1$. The proof of Case $2$ is similar. Suppose that $X$ is a set on $v$ points where $v\equiv 1,3$ (mod $6$). Let $(X,\mathcal{A}_1)$, $(X,\mathcal{A}_2)$, $\cdots$, $(X,\mathcal{A}_t)$ be $t$ Steiner triple systems with $k$ blocks in common. We have two following cases. \\

1) Suppose that each element of $X$ has been appeared in at most one of common blocks. Since $k\geq 2$, suppose that $xyz$ and $uvw$ are two of common blocks. The third element in block containing the pair $xu$ must be distinct in each of $t$ systems and this element  belongs to $X\setminus \{x,y, z, u, v, w \}$. Then $t\leq v-6$.\\

2) Assume that there exists an element $x\in X$ that has been appeared in at least two of common blocks. If the element $x$ belongs to some of uncommon blocks, let $p\in X$ be an element that the pair $xp$   belongs to one of uncommon blocks and $xyz$ and $xuv$ be two of common blocks, then the third element in block containing the pair $xp$ must be distinct in each of $t$ systems and this element  belongs to $X\setminus \{x,p, y, z, u, v\}$. Then $t\leq v-6$. If the element $x$ does not belong to any of uncommon blocks, since $k\neq \frac{v-1}{2}$, let $y$ be an element of $X\setminus \{x\}$ that has been appeared in at least two of common blocks and at least one of uncommon blocks. Suppose that $xyz$ and $yuv$ are the mentioned common blocks and $p$ is an element of $X\setminus \{x,y,z,u,v\}$  such that the pair $yp$ does not belong to any of common blocks and $xpw$  is the common block containing the pair $xp$. Then the third element in block containing the pair $yp$ must be distinct in each of $t$ systems and this element belongs to $X\setminus \{x,y,z,u,v,p,w\}$. Then $t\leq v-7$. \\
\end{proof}
 \end{theorem}
\begin{theorem}\label{7}
If $v\equiv 1,3$ (mod $6$) and $v\geq 13$, then $\mathcal{D}(v,k)\leq v-7$ for $k\geq \lfloor\frac{v}{3}\rfloor +1$ and $k\neq \frac{v-1}{2}$.
\begin{proof}
Suppose that $X$ is a set on $v$ points where $v\equiv 1,3$ (mod $6$). Let $(X,\mathcal{A}_1)$, $(X,\mathcal{A}_2)$, $\cdots$, $(X,\mathcal{A}_t)$ be $t$ Steiner triple systems with $k$ blocks in common. Since $k\geq \lfloor\frac{v}{3}\rfloor +1$, then by the pigeonhole principle there is at least one element  in $ X$ that has been appeared in at least two blocks of $k$ common blocks. The two following cases will be happen.\\ 

1) If there exists an element $x$ that does not belong to any of uncommon blocks, the assertion follows from the proof of Theorem \ref{6}. \\

2) Suppose that each element $x\in X$ belongs to some of uncommon blocks. If there is an element $x$ that has been appeared in at least three of common blocks, obviously $t\leq v-7$.  If for each element $x$, the number of common blocks containing the element $x$ is at most two, let $xuv$ and $xwz$  be two of common blocks. Since $k\geq 5$ and each element has been appeared in at most two of common blocks, then there exists one common block that has at least two elements $y$ and $p$ distinct from the set $\{x,u,v,w,z\}$. Then the third element in block containing the pair $xp$ must be distinct in each of $t$ systems and this element belongs to the set $X\setminus \{x,u,v,w,z,y,p\}$. So $t\leq v-7$.
\end{proof}
 \end{theorem}
\section{Minimum volume and  foundation size of  $\mu$-way $(v,3,2)$ trades}
Hwang in \cite{hwang}  proved the following theorem for $\mu=2$. The proof can be easily obtained by induction on $t$ and using Lemma \ref{32}.
\begin{theorem}\label{ssss}
If $T$ is a $(v,k,t)$ trade, then 
\begin{itemize}
\item[i.]
$|$found(T)$|\geq k+t+1$,
\item[ii.]
the volume of $T$ is at least $2^t$.
\end{itemize}
\end{theorem}
No analogue of  Theorem  \ref{ssss} for $\mu$-way trades ($\mu\geq 3)$ is known. In this section, we determine the minimum volume and minimum foundation size of $\mu$-way $(v,3,2)$ trades for each integer $\mu\geq 3$.    The following theorem help us to obtain the main result.
\begin{theorem}\label{a}
If $T$ is a $\mu$-way $(v,2,1)$ trade, then the volume of $T$ is at least $\lceil\frac{\mu+1}{2}\rceil$.
\begin{proof}
Suppose that $T=(X,\{T_1, T_2, \cdots, T_{\mu}\})$ is a $\mu$-way $(v,2,1)$ trade of volume $m$ and $x_1\in$ found(T). Let $x_{i+1}$ be the element of found(T) such that $x_1x_{i+1}\in T_i$ for $i=1, 2, \cdots, \mu$. The distinct elements $x_3,x_4, \cdots, x_{\mu+1}$ must be appeared in at least one block of  $T_1$ and so at least $\lceil\frac{\mu-1}{2}\rceil$ blocks will be needed. Then $m\geq \lceil\frac{\mu-1}{2}\rceil+1$. On the other hand, the following structure is a $\mu$-way $(v,2,1)$ trade of volume $\lceil\frac{\mu+1}{2}\rceil$ in which $F_1$, $F_2$, $\cdots$, $F_{\mu}$ are the $1$-fractors of complete graph $K_{2\lceil\frac{\mu+1}{2}\rceil}$. Obviously, $|F_i|=\lceil\frac{\mu+1}{2}\rceil$.
\begin{center}
\begin{tabular}{c|c|c|c}
$T_1$&$T_2$& $\cdots$ & $T_{\mu}$\\
\hline
$F_1$&$F_2$& $\cdots$ & $F_{\mu}$\\
\end{tabular}
\end{center}
Then, the minimum volume of a $\mu$-way $(v,2,1)$ trade is $\lceil\frac{\mu+1}{2}\rceil$. Then the proof is complete.
\end{proof}
\end{theorem}
The following corollary immediately follows from Lemma \ref{32} and the above theorem.
\begin{corollary}\label{sssss}
If $T$ is a $\mu$-way $(v,3,2)$ trade, then for each $x\in found(T)$, $r_x\geq \lceil\frac{\mu+1}{2}\rceil$.
\end{corollary}
Now, we want to show that a $\mu$-way $(v,3,2)$ trade  of minimum volume must be a Steiner trade. Before showing this fact, we give  an upper bound for the minimum volume of a $\mu$-way $(v,3,2)$ trade depending on the value of $\mu$.
\begin{lemma}\label{2} Suppose that $\mu=6k+s$ and $\mu\geq 3$  is an integer number  and $T=(X,\{T_1,T_2,\cdots, T_{\mu}\})$  is a $\mu$-way $(v,3,2)$ trade of minimum volume $m$, then
\begin{itemize}
\item[i.]
If  $s=0$ or $1$,  then $m\leq \frac{(\mu+3-s)(\mu+2-s)}{6}$.
\item[ii.]
If $2\leq s\leq 5$ and $\mu\geq 8$, then $m\leq \frac{(\mu+7-s)(\mu+6-s)}{6}$,   $\mu=3$, $m\leq 6$, $\mu=4$, $m\leq 8$ and $\mu=5$, $m\leq 12$.
\end{itemize}
\begin{proof}
For $s=0$ or $1$, by the first part of Theorem \ref{1}, there exist  $\mu+1-s$ disjoint $STS(\mu+3-s)$s that form a $\mu$-way $(v,3,2)$ trade of volume $\frac{(\mu+3-s)(\mu+2-s)}{6}$. Then  the minimum volume of a $\mu$-way $(v,3,2)$ trade for $\mu\equiv 0$ or $1$ (mod $6$)  is at most  $\frac{(\mu+3-s)(\mu+2-s)}{6}$.\\

For $2\leq s\leq 5$ and $\mu\geq 8$,  by the first part of Theorem \ref{1}, there exist $\mu+5-s$ disjoint $STS(\mu+7-s)$s. By choosing $\mu$ of them, we can form a $\mu$-way $(v,3,2)$ trade of volume $\frac{(\mu+7-s)(\mu+6-s)}{6}$. Then for $\mu\equiv 2,3,4,5$  (mod $6$) and $\mu\geq 8$ the minimum volume of a $\mu$-way $(v,3,2)$ trade  is at most $\frac{(\mu+7-s)(\mu+6-s)}{6}$.
For $\mu=3$, by Example \ref{pp} $m\leq 6$, for $\mu=4$, by Example  \ref{pppp} $m\leq 8$, and for $\mu=5$, since $\varphi(9)=7$ then $m\leq 12$.
\end{proof}
\end{lemma}
\begin{theorem}\label{3}
Each $\mu$-way $(v,3,2)$ trade of minimum volume  is a Steiner trade.
\begin{proof}
Suppose that $T=(X,\{T_1,T_2,\cdots, T_{\mu}\})$ is a $\mu$-way $(v,3,2)$ trade of  minimum volume. Let $x$ and $y$ be two elements of found(T) such that $\lambda_{xy}\geq 2$. Suppose that $\alpha_i$ and $\beta_i$ for $i=1,2,\cdots,\mu$ are the elements of found(T) such that $xy\alpha_i$, $xy\beta_i\in T_i$. Then there must exist some blocks in $T_1$ containing the pairs $x\alpha_i$, $x\beta_i$ and $y\alpha_i$, $y\beta_i$ for $i=2,3,\cdots,\mu$. In the best condition, there exist $1$-factors $F_1$ and $F'_1$ on the set $\{\alpha_i,\beta_i: i=2,3,\cdots,\mu\}$ such that $xF_1\subset T_1$ and $yF'_1\subset T_1$. Similarly, there exist $1$-factors $F_k$ and $F'_k$ on the set $\{\alpha_i,\beta_i: i\in \{1,2,\cdots,\mu\}\setminus \{k\}\}$ such that $xF_k\subset T_k$ and $yF_k'\subset T_k$. \\
Suppose that $\mathcal{F}=\bigcup_{i=2}^{\mu}(F_i\cup F'_i)$. Obviously, $|\mathcal{F}\cap (F_1\cup F'_1)|\leq 2(\mu-1)$ and the number of elements of $\mathcal{F}$ by counting the repeated elements is $2(\mu-1)(\mu-1)$. Then
$$|\mathcal{F}\setminus (F_1\cup F'_1)|\geq 2(\mu-1)(\mu-1)-2(\mu-1)=2(\mu-1)(\mu-2)$$
Since each element of $\mathcal{F}$ has been counted at most two times, then the number of  distinct elements of $\mathcal{F}$ disjoint from $F_1$ and $F'_1$ is at least $\frac{2(\mu-1)(\mu-2)}{2}=(\mu-1)(\mu-2)$. Then
$$|T|=|T_1|\geq \frac{(\mu-1)(\mu-2)}{3}+2(\mu-1)+2$$
and this is in contradiction to the result in Lemma \ref{2}.
\end{proof}
\end{theorem}
\begin{theorem}\label{8}
If $T$ is  a $5$-way $(v,3,2)$ trade, then the volume of $T$ is at least $12$.
\begin{proof}
Since there exist seven disjoint STS$(9)$s, then by choosing $5$ of them, we have a $5$-way $(9,3,2)$ trade of volume $12$. Then it is sufficient to show that there does not exist any $5$-way $(v,3,2)$ trade of volume less than $12$. Suppose that   $T=(X,\{T_1,T_2,T_3,T_4,T_5\})$ is a $5$-way $(v,3,2)$ trade of volume  $m<12$. By Theorem \ref{1}, since $\varphi(7)=2$ and $\varphi(8)=4$, then $|$found(T)$|\geq 9$. Since $m<12$, then there must  exist at least one element of found(T) for example $x$ such that $r_x= 3$. Let $F_i$ be the $1$-factor on the set $\{1,2,3,4,5,6\}$ such that $ T_{ix}=xF_i$.\\

Set $\mathcal{B}^i=\bigcup_{j\neq i}F_j$. 
By  $\alpha_i$ and $\beta_i$ we denote the number of blocks in $T'_{ix}$ that contain three elements and one element of $\mathcal{B}^i$, respectively for $i=1,2,3,4,5$.  Since $T$ is a Steiner trade, then
$$\alpha_i\leq 4, \ \ \ 3\alpha_i+\beta_i=12,$$
and since $m\leq 11$, then
$$\alpha_i+\beta_i\leq 8.$$
Therefore $\alpha_i\geq 2$. Now we show that $\alpha_i\neq 3,4$ for each $i$.\\

If for some $i$, $\alpha_i=4$, since there does not exist $5$-way $(v,3,2)$ trade of volume $7$ and $T$ is a Steiner trade, then $T'_{ix}$ contains a block $\{y,w,z\}$ that at most  one element of this block belongs to $\{1,2,\cdots,6\}$. Since $|\{y,w,z\}\cap \{1,2,\cdots,6\}|\leq 1$,  $r_y,r_w,r_z\geq 3$ and $T$ is a Steiner trade, then   $m\geq 3+4+1+2\times 2$ that it is in contradiction to our hypothesis $m\leq 11$. Therefore $\alpha_i\neq 4$.\\

Now assume that for some $i$, $\alpha_i=3$ and $A$ is the set of these three blocks. Since $T$ is a Steiner trade, then nine distinct elements of $\mathcal{B}^i$  have been appeared in  blocks of $A$. By the pigeonhole principle three elements of $\{1,2,\cdots,6\}$ have been appeared in two blocks of $A$ and the remaining  elements in one block of $A$. If $a,b$ and $c$ are three last elements, then there exist three distinct elements $y,w$ and $z$ disjoint from $\{1,2,\cdots,6\}$ such that $\{y,a,b\}$, $\{w,a,c\}$ and $\{z,b,c\}\in T'_{ix}$. Since $r_y,r_w$ and $r_z\geq 3$, then $m\geq 3+3+3+3$ that it is in contradiction to our hypothesis. Therefore $\alpha_i\neq 3$ and so $\alpha_i=2$.\\

Suppose that $A_i$ and $B_i$ are the set of blocks in $T'_{ix}$ that are containing three elements and one element of $\mathcal{B}^i$. Then $|A_i|=2$ and $|B_i|=6$. If   there exists an element $b\in \{1,2,\cdots,6\}$ such that $b$ does not belong to any blocks of $A_i$, then the element $b$ belongs to four blocks of $B_i$ and so $r_b\geq 5$. Since $T$ is a Steiner trade, then  there exist the distinct elements $y_1,y_2,y_3,y_4\notin \{1,2,\cdots,6\}$  such that the pairs $by_i$ for $i=1,2,3,4$ belong to four blocks of $B_1$. Then $|$found(T)$|\geq 11$ and  $m\geq \frac{3\times 10+5}{3}$ that is in contradiction to our hypothesis. \\

Then we conclude if there exists a $5$-way $(v,3,2)$ trade $T$ of volume less then $12$, the blocks in $A_i$ must be in the form $\{a_i,b_i,c_i\}$ and $\{1,2,\cdots,6\}\setminus \{a_i,b_i,c_i\}$, for each $i$.  Then the not appeared pairs of $\mathcal{B}^i$ in $A_i$ form two $1$-factors on the set $\{1,2,3,4,5,6\}$ that we denote them by $F_i$ and $F'_i$ in which $F_i\cap F'_i=\emptyset$. Then there exist two distinct elements $w,z\notin \{1,2,\cdots,6\}$ such that $wF_i,zF'_i\in B_i$. One can easily seen that there  exists a pair $ab$  that has  been appeared in a block of $A_i$ for an unique value  $i$. Suppose that  $xab\in T_{jx}$ ($j\neq i$). Then for each three index $k\neq i,j$, the pair $ab$ must be appeared in a block of $B_k$. By the pigeonhole principle there exist the indices $k$ and $k'$ such that $wab\in B_k$ and $wab\in B_{k'}$ that is a  contradiction. Then the hypothesis $m\leq 11$ is impossible and so $m\geq 12$.
\end{proof}
\end{theorem}
\begin{theorem}If $T$ is a $\mu$-way $(v,3,2)$ trade, then 
\begin{itemize}
\item[1)] for $\mu\equiv 0,4$ (mod $6$) and $\mu\neq 4$ the volume $T$ is at least  $\frac{(\mu+3)(\mu+2)}{6}$ and $|$found(T)$|\geq \mu+3$, for $\mu=4$ is at least $8$ and $|$found(T)$|\geq 8$,
\item[2)] for $\mu\equiv 1,5$ (mod $6$) and $\mu\neq 5$ the volume $T$ is at least $\frac{(\mu+2)(\mu+1)}{6}$ and $|$found(T)$|\geq \mu+2$, for $\mu=5$ is at least $12$ and $|$found(T)$|\geq 9$,
\item[3)] for $\mu\equiv 2$ (mod $6$) the volume $T$ is at least  $\frac{(\mu+4)(\mu+2)}{6}$ and $|$found(T)$|\geq \mu+4$,
\item[4)] for $\mu\equiv 3$ (mod $6$) and $\mu\neq 3$ the volume $T$ is at least   $\frac{(\mu+4)(\mu+3)}{6}$ and $|$found(T)$|\geq \mu+4$,   for $\mu=3$ is at least $6$ and $|$found(T)$|\geq 7$.
\end{itemize}
\begin{proof}
Let $x$ be an element of found(T) with $r_x<m$. Then by Corollary \ref{sssss}, $r_{x}\geq \lceil \frac{\mu+1}{2}\rceil$. By Theorem \ref{3}, Since a $\mu$-way $(v,3,2)$ trade of minimum volume is a Steiner trade, then for all $w,z\in$found(T), $\lambda_{wz}=0$ or $1$.\\

$1)$  Since $\mu$ is even, then $r_x\geq \frac{\mu+2}{2}$. Also $|$found(T)$|\geq 2\lceil \frac{\mu+1}{2}\rceil+1=\mu+3$. Now, we have two cases:\\
If for all $w,z\in $found(T), $\lambda_{wz}=1$, then
\begin{eqnarray}\label{aaaaaa}
m\geq \frac{
\left( {\begin{array}{*{20}c}
\mu+3 \\ 
2 \\ 
\end{array}} \right)
}{3}=\frac{(\mu+3)(\mu+2)}{6}
\end{eqnarray}
 If there exist two elements $w,z\in$found(T) such that $\lambda_{wz}=0$, then $|$found(T)$|\geq 2\lceil \frac{\mu+1}{2}\rceil+1+1=\mu+4$. From the relation $3m=\sum_{x_i\in found(T)}r_{x_i}\geq (\mu+4)\frac{\mu+2}{2}$ we have
 \begin{eqnarray}\label{bbbbb}
  m\geq \frac{(\mu+4)(\mu+2)}{6}.
 \end{eqnarray}
Since by the first part of Theorem \ref{1}, except for $\mu=4$ there exist $\mu+1$ disjoint Steiner triple systems on $\mu+3$ points with $\frac{(\mu+3)(\mu+2)}{6}$ blocks, then the assertion follows from \ref{aaaaaa} and \ref{bbbbb}. For $\mu=4$, since by Theorem \ref{1}, $\varphi(7)=2$, then $|$found(T)$|\geq 8$ and so $m\geq \frac{8\times 3}{3}$. Now, by Example \ref{pppp}, we conclude that $m=8$. \\

 $2)$ In this case for all $x\in$found(T), $r_x\geq \frac{\mu+1}{2}$ and $|$found(T)$|\geq 2\lceil \frac{\mu+1}{2}\rceil+1=\mu+2$. There are two following cases: \\
 If for all $w,z\in $found(T), $\lambda_{wz}=1$, then
 \begin{eqnarray}\label{ccccc}
m\geq \frac{
\left( {\begin{array}{*{20}c}
\mu+2 \\
2 \\
\end{array}} \right)
}{3}=\frac{(\mu+2)(\mu+1)}{6}
\end{eqnarray}
If there exist two elements $w,z\in$found(T) such that $\lambda_{wz}=0$, then $|$found(T)$|\geq 2\lceil \frac{\mu+1}{2}\rceil+1+1=\mu+3$. From the relation $3m=\sum_{x_i\in found(T)}r_{x_i}\geq (\mu+3)\frac{\mu+1}{2}$ we have
\begin{eqnarray}\label{ddddd}
 m\geq \frac{(\mu+3)(\mu+1)}{6}.
\end{eqnarray}
Since by the first part of Theorem \ref{1}, except for $\mu=5$ there exist $\mu$ disjoint Steiner triple systems on $\mu+2$ points with $\frac{(\mu+2)(\mu+1)}{6}$ blocks, then the assertion follows from \ref{ccccc} and \ref{ddddd}. For $\mu=5$,  the claim has beeen proved in Theorem \ref{8}.\\

 $3)$ In this case for all $x\in$found(T), $r_x\geq \frac{\mu+2}{2}$ and $|$found(T)$|\geq 2\lceil \frac{\mu+1}{2}\rceil+1=\mu+3$. Since $\mu+3\equiv 5$ (mod $6$), then there must be  $w,z\in$found(T) such that $\lambda_{wz}=0$. Then $|$found(T)$|\geq 2\lceil \frac{\mu+1}{2}\rceil+1+1=\mu+4$. From the relation $3m=\sum_{x_i\in found(T)}r_{x_i}\geq (\mu+4)\frac{\mu+2}{2}$ we have
 \begin{eqnarray}\label{eeeee}
  m\geq \frac{(\mu+4)(\mu+2)}{6}.
 \end{eqnarray}
 Since $\mu+4\equiv 0$ (mod $6$), by the second part of Theorem \ref{1}, there exist $\mu$ disjoint mutually balanced $(2,3)$-packings  on $\mu+4$ points with $\frac{(\mu+4)(\mu+2)}{6}$ blocks, then the assertion follows from \ref{eeeee}.\\


 $4)$ In this case, for all $x\in$ found(T), $r_x\geq \frac{\mu+1}{2}$ and $|$found(T)$|\geq \lceil \frac{\mu+1}{2}\rceil+1=\mu+2$. By Theorem \ref{0}, since $\varphi(\mu+2)\leq \mu-2$ and $\varphi(\mu+3)\leq \mu-1$, then $|$found(T)$|\geq \mu+4$. Since $\mu+4\equiv 1$ (mod $6$), by the first part of Theorem \ref{1} except for $\mu=3$ there exist $\mu+2$ disjoint Steiner triple systems on $\mu+4$ points and with $\frac{(\mu+4)(\mu+3)}{6}$ blocks. Then
 $$m\leq \frac{(\mu+4)(\mu+3)}{6}.$$
 On the other hand since for all $x\in$found(T), $r_x\geq \frac{\mu+1}{2}$,
 $$\frac{(\mu+4)(\mu+3)}{2}\geq 3m=\sum_{x\in found(T)}r_x\geq \frac{\mu+1}{2}|found(T)|,$$
 then $|$found(T)$|\leq \mu+6$.\\

 For the case $|$found(T)$|=\mu+4$, since $\mu+4\equiv 1$ (mod $6$), by Theorems \ref{4} and \ref{6} $\mathcal{D}(\mu+4,k)\leq \mu+4-5=\mu-1$ for $k\geq 1$, then the minimum volume of $\mu$-way $(\mu+4,3,2)$ trade in this case is the number of blocks in a  STS($\mu+4$), i.e. $\frac{(\mu+4)(\mu+3)}{6}$.\\

 In the case $|$found(T)$|=\mu+5$, since $\mu+5\equiv 2$ (mod $6$), by Theorems \ref{4} and \ref{6} $\mathcal{D}(\mu+5,k)\leq \mu+5-6=\mu-1$ for $k\geq 1$ and $\varphi(\mu+5)=\mu+1$, then the minimum volume of $\mu$-way $(\mu+5,3,2)$ trade in this case is the number of blocks in a  maximum packings on $\mu+5$ points, i.e. $\frac{(\mu+5)(\mu+3)}{6}$.\\

 For the case $|$found(T)$|=\mu+6$, Since by Theorem  \ref{7}, $\mathcal{D}(\mu+6,k)\leq \mu+6-7=\mu-1$ for $k\geq \frac{\mu+6}{3}+1$, then the minimum volume of $T$ is at least
 $$\frac{(\mu+6)(\mu+5)}{6}-\frac{\mu+6}{3}=\frac{(\mu+6)(\mu+3)}{6}.$$
 By comparison the three numbers we conclude the minimum volume of a $\mu$-way $(v,3,2)$ trade for $\mu\equiv 3$ (mod $6$) and $\mu\neq 3$, is  $\frac{(\mu+4)(\mu+3)}{6}.$ For $\mu=3$, the claim has been proved in \cite{rashidi}.

\end{proof}
\end{theorem}


\begin{thebibliography}{10}
\bibitem{amjadi}
 H. Amjadi, N. Soltankhah, \emph{On the existence of $d$-homogeneous $3$-way
Steiner trades}, Utilitas. Math. (to appear).
 \bibitem{asgari}
 M. Asgary, N. Soltankhah, \emph{On the non-existence of some Steiner $t-(v,k)$ trades of certain volumes}, Utilitas. Math.  \textbf{79} (2009) 277-283.
  \bibitem{cao1}
  H. Cao, L. Ji, L. Zhu, \emph{Large sets of disjoint packings on $6k+5$ points}, J. Combin. Theory (A). \textbf{108} (2004) 169-183.
  \bibitem{cao2}
  H. Cao, J. Lei, L. Zhu, \emph{Large sets of disjoint group-divisible designs with block size three and type $2^n4^1$}, J. Combin. Des. \textbf{9} (2001) 285-296.
  \bibitem{cao3}
  H. Cao, J. Lei, L. Zhu, \emph{Further results on Large sets of disjoint group-divisible designs with block size three and type $2^n4^1$}, J. Combin. Des. \textbf{11} (2003) 24-35.
  \bibitem{cao4}
  H. Cao, J. Lei, L. Zhu, \emph{Constructions of Large sets of disjoint group-divisible designs  LS($2^n4^1$) using a generalization of LS($2^n$)}, Discrete Math. \textbf{338} (2015) 1449-1559.
  \bibitem{chang}
  Y. Chang, L. Ji, H. Zheng, \emph{A completion of   LS($2^n4^1$)}, Discrete Math. (to appear).
  \bibitem{chen}
  D. Chen, C. C. Lindner, D. R. Stinson, \emph{Further results on Large sets of disjoint group-divisible designs}, Discrete Math.  \textbf{110} (1992) 35-42.
  \bibitem{chen2}
  D. Chen, R. G. Stanton, D. R. Stinson, \emph{Disjoint packings on $6k+5$ points}, Utilitas. Math. \textbf{40} (1991) 129-138.
  \bibitem{chen3}
  D. Chen, D.R. Stinson, \emph{Recent results on combinatorial constructions for threshold schemes}, Australas. J. Combin.  \textbf{1} (1990) 29-48.
  \bibitem{chen4}
  D. Chen, D. R. Stinson, \emph{On the construction of Large sets of disjoint group-divisible designs}, Ars Combin. \textbf{35} (1993) 103-115. 
  \bibitem{forbes}
  A. D. Forbes, M. J. Grannell, T. S. Griggs, \emph{Configurations and trades in Steiner triple systems}, Australas. J. Combin.  \textbf{29} (2004) 75-84.
  \bibitem{hedayat}
 A. S. Hedayat, \emph{The theory of trade-off for $t$-designs}, In: D. Ray-Chaudhuri:
Coding theory and design theory, Part II: Design Theory, IMA Vol. Math.
Appl. \textbf{21} (1990) 101-126.
\bibitem{hedayat2}
  A. S. Hedayat, B. Khosrovshahi, \emph{Trades, Handbook of combinatorial designs}, CRC Press.  (2007) 644-648.
  \bibitem{hwang}
H.~L. Hwang, \emph{On the structure of {$(v,k,t)$} trades}, J. Statist. Plann.
  Inference. \textbf{13} (1986), no.~2, 179-191.
  \bibitem{ji}
L. Ji, \emph{Existence of large sets disjoint group divisible designs with block size three and type $2^n4^1$}, J. Combin. Des. \textbf{13} (2005), 302-312.
  \bibitem{lei1}
  J. Lei, \emph{Completing the spectrum for LGDD($m^v$)}, J. Combin. Des. \textbf{5} (1997) 1-11.
   \bibitem{lindner}
  C. C. Lindner, C. A. Rodger, \emph{Maximum packings and minimum coverings, Design theory}, CRC Press. (1997) 53-63.
  \bibitem{lu1}
  J. X. Lu, \emph{On large sets of disjoint Steiner triple systems I, II, and III}, J. Combin. Theory (A). \textbf{34} (1983) 140-146, 147-155, 156-182.
  \bibitem{lu2}
  J. X. Lu, \emph{On large sets of disjoint Steiner triple systems IV, V, and VI}, J. Combin. Theory (A). \textbf{37} (1984) 136-163, 164-188, 189-192.
  \bibitem{mahmood}
  E. S. Mahmoodian, N. Soltankhah, \emph{On the existence of $(v,k,t)$ trades}, Australas. J. Combin. \textbf{6} (1992) 279-291.
  \bibitem{rashidi}
  S. Rashidi, N. Soltankhah, \emph{On the possible volume of $\mu-(v,k,t)$ trades}, Bull. Iranian Math. Soc. \textbf{40} (2014) 1387-1401.
  \bibitem{schellenberg}
  P. J. Schellenberg, D. R. Stinson, \emph{Threshold schemes from combinatorial designs}, J. Combin. Math. Combin. Comput \textbf{5} (1989) 143-160.
  \bibitem{teirlink}
  L. Teirlinck,  \emph{A completion of Lu's determination of the spectrum of large sets of disjoint Steiner triple systems}, J. Combin. Theory (A). \textbf{57} (1991) 302-305.
  \end{thebibliography}
\end{document}